\documentclass{article}

\usepackage{amsmath,amsthm}

\newtheorem{theorem}{Theorem}

\newtheorem{proposition}{Proposition}

\theoremstyle{remark}

\theoremstyle{definition}

\theoremstyle{definition}

\begin{document}

\title{DYNAMICS OF CONTROLLED HYBRID SYSTEMS
OF AERIAL CABLE-WAYS\footnote{Accepted (15-May-2006)
to the Proceedings of the
{\it International Conference of Hybrid Systems and Applications},
The University of Louisiana, Lafayette, LA, USA, May 22-26 2006,
to be published in the journal \emph{Nonlinear Analysis:
Hybrid Systems and Applications}.}}

\author{Olena V. Mul${}^{1)}$\thanks{Supported by FCT
(the Portuguese Foundation for Science
and Technology), fellowship SFRH/BPD/14946/2004.}
        \and
        Delfim F. M. Torres${}^{1)}$
        \and
        Volodymyr P. Kravchenko${}^{2)}$
        }

\date{${}^{1)}$ Department of Mathematics, University of Aveiro\\
      Campus Universit\'{a}rio de Santiago\\
      3810-193 Aveiro, Portugal\\
      Email: \texttt{\{olena,\,delfim\}@mat.ua.pt}\\
      [0.3cm]
      ${}^{2)}$ Physico-Technological Institute of Metals and Alloys\\
      of the National Academy of Sciences of Ukraine\\
      Academician Vernadsky Avenue, 34/1\\
      03680 Kiev-142, Ukraine\\
      Email: \texttt{v$_{-}$kram@i.com.ua}}

\maketitle


\begin{abstract}
Dynamics of the hybrid systems of aerial cable-ways is investigated.
The eigenvalue problems are considered for such hybrid systems with
different assumptions. An overview of different methods for
eigenvalue problems is given. In the research, the method of the
normal fundamental systems is applied, which turns out to be very
effective for the considered problems. Changes of dynamical
characteristics of the systems depending on the controlled parameter
are studied.
\end{abstract}


\smallskip


\smallskip

\smallskip



\section{Introduction}

Non-stationary boundary problems with continuous-discrete parameters
in the theory of vibrations belong to the wide class of
continuous-discrete boundary problems, whose parameters change with
time. These problems are of great importance both theoretically and
in applications \cite{Kach-teor}. Typically, their study combines
solving non-stationary problems with continuous parameters
\cite{myNY98,myDresd,Sloss} and non-stationary problems with
continuous-discrete parameters \cite{Kach-teor,my_nonlin_anal}.

The dynamics for many elastic systems is described by partial
differential equations with non-stationary continuous-discrete
coefficients. Moreover, the solutions of those equations have to
satisfy not only some boundary conditions but also some conjugation
conditions, which considerably complicate the analytical or
numerical study. The use of normal fundamental functions
\cite{Norm-fund} allows to find recurrence formulas, by means of
which the conjugation conditions are taken into account. This
simplifies considerably the resolution of the non-stationary
continuous-discrete boundary problems.

In this paper dynamics of some controlled systems of aerial
cable-ways with continuous and discrete parameters is studied. For
research, mathematical models in the form of the eigenvalue problems
are used. In \S\ref{sec:OMEP} for similar eigenvalue problems it is
given an overview of different known methods, such as variational
methods, the Raleigh's energy method and methods of integral
equations. In \S\ref{sec:NFFM} an application of the method of the
normal fundamental systems, which turns out to be the most effective
one for the hybrid eigenvalue problems, is described. In
\S\ref{sec:EPAC} we consider several eigenvalue problems that model
the dynamics of aerial cable-ways at different assumptions about the
system which is: only the cable without discrete loads, the cable
with fixed loads, the cable with loads moving along it at a steady
speed, or the cable with fixed discrete loads moving at a steady
speed. The method of the normal fundamental systems is applied to
these problems, and the dependence of the dynamical characteristics
of the systems on the controlled parameter, which is the speed of
the system motion, is studied.


\section{An Overview of Methods for Eigenvalue Problems}
\label{sec:OMEP}

The mathematical problems on eigenvalues are widespread both in theoretical
and applied mathematics. They serve as models
for many concrete problems in physics, engineering and biology.
Most of these applied problems may be reduced to the study
of boundary problems for ordinary or partial differential equations,
where almost always the main question is reduced to the
determination of the eigenvalues and eigenfunctions \cite{Akulenko}.

There exist many different non-stationary boundary problems with
continuous and discrete parameters which are mathematical models for
a wide class of controlled dynamical systems \cite{Kach-teor}.
Dynamical characteristics of such hybrid systems may be controlled
by changing the system parameters within certain limits. Vibrations
in such systems are described both by partial differential equations
with variable coefficients and by homogeneous Fredholm integral
equation, with a Stieltjes type of integral.

In order to determine the eigenvalues of the vibration boundary problems,
different variational methods, as well as methods of
differential and integral equations, may be used \cite{Gould,Motreanu}.

An eigenvalue problem may be reduced to a variational problem, viz,
to the search of the minimum of a certain functional
\cite{Norm-fund}. One of the most well-known variational method is
the Ritz method, which can be applied to the problem
\begin{equation}
\label{probl-min} \min_{u} (Au,u)
\end{equation}
subject to the additional condition
\begin{equation}
\label{probl-min2} (u,u)=1 \, ,
\end{equation}
where $A$ is a positive operator and $u$ is a continuous function.
In order to solve problem \eqref{probl-min}-\eqref{probl-min2} by
the Ritz method, one begins by choosing some sequence of basic
functions $\varphi_{n}$, $n=1,\, 2,\, ...$. They have to satisfy the
following three conditions: all functions $\varphi_{n}$ belong to
the operator definition domain, they are linearly independent for
any $n$, and their system $\{\varphi_{n}\}$ is complete. Then, one
searches the approximate solution in the form
\begin{equation}
\label{appr-sol} u_{n}=\sum\limits_{k=1}^n a_{k} \varphi_{k}\, ,
\end{equation}
where the $a_{k}$'s are unknown coefficients, with
$k=\overline{1,n}$. The coefficients $a_{k}$'s should be selected in
a way to satisfy condition \eqref{probl-min2} for the solution
$u_{n}$ and to make the value $(Au_{n},u_{n})$ minimal. The Ritz
method can be applied only to continuous problems, and is not
applicable to the problems with continuous and discrete distribution
of masses we are interested in.

There exist many other variational methods, for example, the
Bubnov-Galerkin method, which is a generalization of the Ritz
approach for operators that are not obligatory positive. However, an
application of the Bubnov-Galerkin method for continuous-discrete
problems results in significant difficulties on the choice of the
coordinate functions for the complex domains. Other variational
methods are also inappropriate to the study of continuous-discrete
problems, because of the complexity of their application. Another
disadvantage of all the variational methods is that the obtained
approximations for the eigenvalues always exceed their real values.

The Raleigh's energy method is also usually applied to solve
eigenvalue problems. The idea of this approach is as follows. If the
mode of vibrations is known, we can always find the frequency of
free vibrations from the constancy of the sum of the kinetic and
potential energy. Then, with the known frequency, it is easy to find
the eigenvalues. However, the problem is to choose the mode of
vibrations. Raleigh proposed to take the shape of the system static
deformation from the applied load as the first mode of vibrations.
Since the selection of the next vibration mode is complicated, the
Raleigh's method is convenient only for the determination of the
first frequency and the first eigenvalue. This is not enough for our
purposes.

Eigenvalue problems may be also reduced to the determination of
eigenvalues of integral equations of the form
\begin{equation}
\label{int-eq} y(x)=\lambda \int_{a}^{b}r(x)G(x,s)y(s)ds \, ,
\end{equation}
where $G(x,s)$ is the Green function; $r(x)$ is a continuous
positive function; $y(x)$ is an unknown function; $\lambda$ is an
unknown eigenvalue; $a$ and $b$ are known constants.

Equation \eqref{int-eq} can be transformed to an integral equation
with a symmetric kernel. The advantage of the eigenvalue problem
method, which deals with the determination of eigenvalues and
eigenfunctions in the case of a symmetric kernel, is that the same
integral equations describe different vibration processes, such as
flexural, torsional, longitudinal or flexural-and-torsional ones.
Besides, in spite of the different physical sense for the Green's
function, the integral equations are the same for the systems with
continuous and with continuous-discrete parameters. On the base of
the integral methods it is possible to create an universal technique
for solving eigenvalue problems for discrete and discrete-continuous
systems. The big problem with the integral methods, which explains
why they are not usually used in applications, is the great
complexity in the construction of the Green's function.

Thus, most of the known methods for eigenvalue problems cannot be
applied in an effective way to the case of continuous-discrete
systems. In the next section we consider an efficient numerical
method for such hybrid eigenvalue problems.


\section{The Normal Fundamental Function Method}
\label{sec:NFFM}

Here we briefly illustrate the normal fundamental function method by
using the following system of partial differential equations with variable
coefficients \cite{my_Japan}:
\begin{equation}
\label{eq:1}
\begin{gathered}
\frac{{\partial z_k }}{{\partial y}} = \sum\limits_{j=1}^N A_{kj}
(y)z_j +  \sum\limits_{j = 1}^N B_{kj}(y)\frac{{\partial ^2 z_j
}}{{\partial t^2 }}
+  \sum\limits_{j = 1}^N C_{kj} (y)\frac{{\partial z_j }}{{\partial t}} \, , \\
y_0  \le y \le y_n \, , \quad k = \overline{1,N} \, ,
\end{gathered}
\end{equation}
satisfying the linear homogeneous boundary conditions
\begin{equation}
\label{eq:3}
\begin{gathered}
z_r(0) = 0 \, , \quad r = \overline{1,m} \, , \\
z_s (l) = 0 \, , \quad s = \overline{m + 1,N} \, , \quad N = 2m,
\end{gathered}
\end{equation}
and the conjugation conditions
\begin{multline}
\label{eq:2} \sum\limits_{k = 1}^N a_{k,p}^{(i + 1)} (y_i )z_k^{(i +
1)} (y_i ) +  \sum\limits_{k = 1}^N {b_{k,p}^{(i + 1)} (y_i)
\frac{{\partial ^2 z_k^{(i + 1)} (y_i )}}{{\partial t^2 }}} \\
+ \sum\limits_{k = 1}^N d_{k,p}^{(i + 1)}(y_i )\frac{{\partial
z_k^{(i + 1)} (y_i )}}{{\partial t}}
=  \sum\limits_{k = 1}^N {a_{k,p}^{(i)} (y_i )z_k^{(i)} (y_i )} \\
+ \sum\limits_{k = 1}^N {b_{k,p}^{(i)} (y_i) \frac{{\partial ^2
z_k^{(i)}(y_i )}}{{\partial t^2 }}} + \sum\limits_{k = 1}^N
{d_{k,p}^{(i)}(y_i )\frac{{\partial z_k^{(i)} (y_i )}}{{\partial
t}}} \, ,
\end{multline}
where $z_j(y,t)$ is the deviation function of the vibration system
from an equilibrium position; $A_{kj}(y)$, $B_{kj}(y)$ and
$C_{kj}(y)$ are real-valued piecewise continuous functions bounded
on the interval $[y_0 ,y_n ]$ for $k,j = \overline{1,N}$; $y$ is the
coordinate of the system point; $t$ is the time variable; $y = y_i$
are the discontinuity points of functions $A_{kj}(y)$, $B_{kj}(y)$
and $C_{kj}(y)$ for $i = \overline{1,n - 1}$; $p = \overline{1,N}$;
the superscript $i$ means the value of the corresponding function on
the $i$-th interval.

The continuous-discrete boundary problem \eqref{eq:1}--\eqref{eq:2}
describes a large number of different dynamical systems
\cite{Norm-fund,Kach-teor}.

For integration of system \eqref{eq:1}--\eqref{eq:2} we apply the
Fourier method of separation of variables, searching the solution
in the form
\begin{equation*}
\label{eq:4} z_k (y,t) = \varphi_k (y)e^{\lambda t} \, .
\end{equation*}

We are interested in finding the values for the parameter $\lambda$,
for which there exist non-trivial solutions $\varphi_k (y)$ of the
system \eqref{eq:1}-\eqref{eq:2}. Such values $\lambda$ are
eigenvalues of the boundary problem, and the corresponding solutions
$\varphi_k (y)$ are eigenfunctions. The eigenvalues for the
vibration system, described by equations \eqref{eq:1}-\eqref{eq:2},
are complex quantities $\lambda  = q + ip$, where $i=\sqrt{-1}$. The
vibrations have increasing or decreasing amplitude, depending on the
sign of the real part $q$ of the eigenvalue.

In many practical vibration theory problems it is important to know
eigenvalues at which vibrations with a constant amplitude occur.
That is only possible at an imaginary eigenvalue $\lambda = ip$,
where $p$ is the eigenfrequency of the vibrations of the considered
system.

Thus, we assume that the system makes harmonic vibrations with
constant amplitude and the form of solution is
\begin{equation}
\label{eq:5} z_k(y,t) = [\varphi _k^{(1)} (y) + \varphi _k^{(2)}
(y)]e^{ipt} \, ,
\end{equation}
where $\varphi_k^{(1)} (y)$ and $\varphi _k^{(2)}(y)$ are
real-valued functions.
Then the problem \eqref{eq:1}-\eqref{eq:2} is reduced to the normal
system of linear ordinary differential equations
\begin{equation}
\label{eq:6} \frac{{d\varphi_k }}{{dy}} = \sum\limits_{j = 1}^N
{a_{kj} (y)\varphi_j} \, , \quad k = \overline{1,N} \, ,
\end{equation}
satisfying $N$ linear homogeneous boundary conditions and also
linear conjugation conditions at the discontinuity points,
where, for simplicity, $\varphi_{k}$ includes both $\varphi _k^{(1)}$
at the values $k=\overline{1,m}$ and $\varphi _k^{(2)}$ at the values
$k=\overline{m+1,N}$;
$a_{kj}(y)$ are bounded piecewise-continuous functions
on the interval $[y_0 ,y_n]$; $k,j = \overline{1,N}$.

Now, for each interval $y_{i - 1} \le y \le y_i$,
$i = \overline{1,n}$, we can solve a Cauchy problem for the system
\eqref{eq:6}  with the initial conditions at points $y = y_{i - 1}$
as follows:
\begin{equation}
\label{eq:10} \Phi ^{(i)} (y_{i - 1} ) = E \, , \quad i =
\overline{1,n} \, ,
\end{equation}
where $E$ is the unit matrix.
According to Picard's Theorem, each such problem has a unique
solution. We can apply some known numerical method, for
example the Runge-Kutta method, to find $N$ linearly independent
solutions of the system \eqref{eq:6}. Such fundamental system of
solutions
\begin{equation}
\label{eq:9} \Phi ^{(i)} (y) = \left| {\begin{array}{*{20}c}
   {\varphi _{1,1}^{(i)} (y)} & {...} & {\varphi _{1,N}^{(i)} (y)}  \\
   {...} & {...} & {...}  \\
   {\varphi _{N,1}^{(i)} (y)} & {...} & {\varphi _{N,N}^{(i)} (y)}  \\
\end{array}} \right| \, ,
\quad y_{i - 1}  \le y \le y_i \, ,
 \quad i = \overline{1,n} \, ,
\end{equation}
with initial conditions \eqref{eq:10} is called a normal fundamental
system \cite{Norm-fund}. Here each function $\varphi_{k,j}^{(i)}(y)$
is defined and continuous at $[y_{i - 1} ,y_i]$; $k$ is a solution
number; $j$ is a function number.
The general solution of the system \eqref{eq:6} at the $i$-th
interval may be written with the help of the normal fundamental
system of solutions as
\begin{equation*}
\label{eq:11} \varphi_k^{(i)} (y) = \sum\limits_{j = 1}^N {C_j^{(i)}
\varphi_{j,k}^{(i)} (y)} \, , \quad y_{i - 1} \le y \le y_i \, ,
\quad k = \overline{1,N} \, , \quad i = \overline{1,n} \, ,
\end{equation*}
where $C_j^{(i)}$ are unknown constants. Imposing the conjugation
and boundary conditions, we arrive to a homogeneous system of $m$
linear algebraic equations in the coefficients $C_j^{(1)}$, from
which one can write a necessary and sufficient condition for the
existence of a non-trivial solution of the boundary problem
\eqref{eq:1}-\eqref{eq:2} as follows:
\begin{equation}
\label{eq:14} D = \det \left| {\sum\limits_{j = 1}^N {u_{j,q}^{(n)}
\varphi _{j,s}^{(n)} (y_n )} } \right| = 0 \, ,
\end{equation}
where the coefficients $u_{j,q}^{(i)}$ for $i = \overline{1,n}$ are
given by recurrence; $s,q = \overline{m + 1,N}$.
Taking into account the dependencies of the functions $\varphi
_{j,s}^{(n)}(y_n)$ on the value $\lambda$, we determine the
eigenfrequencies of the vibration system as roots of the equation
\eqref{eq:14}.

Given a concrete system, the necessary and sufficient condition
\eqref{eq:14} provides a general method to investigate the frequency
spectrum of possible vibrations and its dependence on the different
parameters of the system. This is the main advantage of the
numerical method presented here: the application of the normal
fundamental functions to an arbitrary but finite number of discrete
characteristics is reduced to the determination of zeros of a
function, which is a determinant; the order of this determinant
depends only on the number of the boundary conditions on the
integration interval and does not depend on the number of the
discrete characteristics in the system; the elements of the
determinant are calculated by given recurrence formulas, which is
very convenient for numerical solving with computers.

In the next section we apply the proposed method to the study of the
eigenvalue problems for some aerial cable-ways.


\section{Eigenvalue Problems for Aerial Cable-Ways}
\label{sec:EPAC}

In this paper we analyze some hybrid controlled systems of aerial
cable-ways, which carry some discrete loads. The problem is to study
the dynamics of such systems as well as changes of their dynamical
characteristics depending on different system parameters and, first
of all, depending on the controlled parameter. In this research, the
method of the normal fundamental systems described in
\S\ref{sec:NFFM} is applied, which gives the best results for this
class of eigenvalue problems.

\subsection{The Basic Problem of Transverse Vibrations of the Cable}

First, let us consider only transverse vibrations of the cable. We
assume that the cable has a uniform linear density. Then we have the
next well-known continuous boundary problem \cite{Kach-teor}:
\begin{equation}
\label{eq:motion-disc1} {\rho\frac{\partial^{2}u}{\partial
t^{2}}=T\frac{\partial^{2}u}{\partial x^{2}} \, ,}
\end{equation}
\begin{equation}
\label{bound-cond-disc1} {u(0,t)=u(l,t)=0 \, ,}
\end{equation}
where $t$ is the time variable; $x$ is the coordinate of some point
of the cable; $\rho$ is the mass per unit length; $T$ is the tension
in the cable; $l$ is the length of the cable; $u(x,t)$ is the
function of deviation from equilibrium position for the point
$x\in[0,l]$ at the time moment $t$.

The classical Fourier technique of separation of variables may be
used for such simplified problem. We assume that the solution of the
equation \eqref{eq:motion-disc1} may be presented as a product of
two functions where the first function depends only on the variable
$x$ and the second one depends only on the variable $t$, i.e.
\begin{equation}
\label{solut-form1} {u(x,t)=X(x)F(t) \, .}
\end{equation}
Then, after separating variables, we will have
\begin{equation}
\label{ravenstvo} {\frac{F''(t)}{a^{2}F(t)} = \frac{X''(x)}{X(x)} \,
,}
\end{equation}
where $\displaystyle a=\sqrt{\frac{T}{\rho}}$. The left-hand side of
the relation \eqref{ravenstvo} depends only on $t$ and the
right-hand side depends only on $x$. Consequently, each side of this
equality is equal to some constant $C$, and the following equations
can be written:

\begin{equation}
\label{equation-form-t} {F''(t)- a^{2}CF(t)=0 \, ,}
\end{equation}

\begin{equation}
\label{equation-form-x} {X''(x)- CX(x)=0 \, .}
\end{equation}
We are interested in finding non-trivial solutions of the equation
\eqref{equation-form-x} satisfying boundary conditions
\begin{equation}
\label{conditions-form-x} {X(0)=X(l)=0 \, ,}
\end{equation}
which obviously follow from the conditions \eqref{bound-cond-disc1}.

It is well known \cite{Norm-fund} that for $C \geq 0$ the equation
\eqref{equation-form-x} has only the trivial solution $X(x)\equiv
0$; and that for $C < 0$ non-trivial solutions exist. Thus, with the
notation $C=-\lambda^{2}$, we can write the equation
\eqref{equation-form-x} in the form

\begin{equation}
\label{eqn-form-x1} {X''(x)+ \lambda^{2}X(x)=0 \, .}
\end{equation}
The values $\lambda^{2}$, at which the initial boundary problem
\eqref{eq:motion-disc1}--\eqref{bound-cond-disc1} has non-trivial
solutions, are called the eigenvalues of this boundary problem.
It is easy to show that the eigenvalues may be found from the
relationship
\begin{equation}
\label{eq:determ-frq2} \frac{\sin(\lambda l)}{\lambda} = 0 \, .
\end{equation}
The general solution of the equation \eqref{eqn-form-x1}, with the
accuracy to some constant, is $X(x)=sin\left(\frac{\pi k \,
x}{l}\right)$ (where $k=1, 2, \ldots$ from here forwards), and it
takes place at the values $\lambda^{2}=\lambda^{2}_{k}=
\frac{\pi^{2}k^{2}}{l^{2}}$, which consequently are the eigenvalues
of the boundary problem
\eqref{eq:motion-disc1}--\eqref{bound-cond-disc1}.

Thus, we can see that in our case the eigenvalues $\lambda^{2}_{k}$
are positive and there exists an infinite set of them. Besides,
$\lambda^{2}_{1}<\lambda^{2}_{2}<\cdots<\lambda^{2}_{k}<\cdots$ and
$\lambda^{2}_{\infty} \rightarrow \infty$. It is found also that
$\lambda^{2}_{k}\neq 0$. These important characteristics of the
eigenvalues depend on the kind of differential equations, which
describe vibrations, and on the boundary conditions. It is necessary
to point out that in some cases, for some boundary problems, it is
possible to determine qualitative characteristics of the eigenvalues
without their numerical calculation.

A particular solution of the equation \eqref{eq:motion-disc1} may
be written as follows:

\begin{equation}
\label{part-solut} u_{k}(x,t) = X_{k}(x)F_{k}(t)=sin \left(
\frac{\pi k
x}{l}\right)\left(\alpha_{k}cos\left(a\lambda_{k}t\right)
+\beta_{k}sin\left(a\lambda_{k}t\right)\right) \, ,
\end{equation}
where $\alpha_{k}$ and $\beta_{k}$ are unknown coefficients which
have to be determined from the initial conditions of the concrete
boundary problem. It is easy to see, from the particular solution
\eqref{part-solut}, that the natural frequencies of vibrations of
the considered cable are given by $\omega_{k}=a
\lambda_{k}=\frac{\pi k \, a}{l}$, $k=1, 2, \ldots$

Thus, if we know the eigenvalues of the boundary problem
\eqref{eq:motion-disc1}-\eqref{bound-cond-disc1},
we can always determine its natural frequencies, which are very important
physical characteristics of the system under consideration;
vice versa, if we know the natural frequencies, we can
find the eigenvalues of the problem.

\subsection{The Problem for the Case of Fixed Loads}

Now let us consider the mathematical model of the cable with
fixed loads. In this case we have the boundary problem
\begin{equation}
\label{eq:motn-disc1} {\rho\frac{\partial^{2}u_{i}}{\partial
t^{2}}=T\frac{\partial^{2}u_{i}}{\partial x^{2}} \, ,}
\end{equation}
\begin{equation}
\label{bound-cond-1} {u_{1}(0,t)=u_{n}(l,t)=0 \, ,}
\end{equation}
where $l_{i}$ is the coordinate of the arbitrary point at which
the discrete load with the mass $m_{i}$ is fixed,
$i=\overline{1,n-1}$, $l_{0}=0$, $l_{n}=l$; $u_{i}(x,t)$ are the
functions of deviation from equilibrium position for the point
$x\in[l_{i-1},l_{i}]$ at the moment $t$, $i=\overline{1,n}$. The
solution of the problem
\eqref{eq:motn-disc1}--\eqref{bound-cond-1} must satisfy the
following conjugation conditions:
\begin{equation}
\label{conjug-cond1-disc} {u_{i+1}(l_{i},t)=u_{i}(l_{i},t)} \, ,
\end{equation}
\begin{equation}
\label{conjug-cond2-disc} {m_{i}\frac{\partial^{2}u_{i}}{\partial
t^{2}}|_{x=l_{i}}=T\left.\left(\frac{\partial u_{i+1}}{\partial
x}-\frac{\partial u_{i}}{\partial x}\right)\right|_{x=l_{i}} \, ,}
\end{equation}
where $m_{i}$ are the masses of the discrete loads at arbitrary
points $x=l_{i}$ for all $i=\overline{1,n-1}$.

Let us use the numerical method of the normal fundamental system of
solutions. An application of this method demands the presentation of
the problem as a system of ordinary differential equations of the
first order in normal form, satisfying some boundary conditions. In
order to reduce \eqref{eqn-form-x1} to the necessary form, we
introduce new functions
\begin{equation}
\label{eq:new-funct} T_{i}(x) = X_i(x) \, , \quad S_{i}(x) = X'_i(x)
\, ,
\end{equation}
and write the desired system of differential equations in normal
form as
\begin{equation}
\label{eq:norm-syst}
\begin{cases}
T'_{i} = S_{i} \, , \\
S'_{i} = -\lambda^{2}T_{i} \, ,
\end{cases}
\end{equation}
with the boundary conditions as follows:
\begin{equation}
\label{eq:bc1:gamma} \bar{x} = 0 : \quad T_{1} = 0 \, ,
\end{equation}
\begin{equation}
\label{eq:bc2:gamma}
\begin{split}
\bar{x} = l : \quad T_{n} = 0 \, .
\end{split}
\end{equation}
Using the functions of the normal fundamental system of solutions, we
can present the general solution of the equations
\eqref{eqn-form-x1} at the intervals of continuity in the form
\begin{equation}
\label{eqn-form-x} {X_{i}(x)=
A_{i}S_{i}(x-l_{i-1})+B_{i}T_{i}(x-l_{i-1}) \, , }
\end{equation}
where
\begin{equation*}
S_{i}(x-l_{i-1})=cos(\lambda(x-l_{i-1})) \, , \quad
T_{i}(x-l_{i-1})=\frac{sin(\lambda(x-l_{i-1}))}{\lambda} \, ,
\end{equation*}
and $S''_{i}=-\lambda^{2}S_{i}$, $T''_{i}=-\lambda^{2}T_{i}$.

In the case of only one discrete load,
we obtain the next equation for determination of the
eigenvalues:
\begin{equation}
\label{eq:determ-freq2} {\frac{\sin(\lambda
l)}{\lambda}-\frac{m_{1}}{\rho} \sin (\lambda b_{1}) \sin (\lambda
b_{2})=0 \, .}
\end{equation}
Comparing the formulas \eqref{eq:determ-freq2} and
\eqref{eq:determ-frq2}, it is easy to see that the presence of the
load at the object decreases the corresponding eigenvalues.
If the load is located in the middle of the integration interval,
then the equation for determination of the eigenvalues is reduced to

\begin{equation}
\label{eq:determ-freq-m} {\frac{\sin(\lambda
l)}{\lambda}-\frac{m_{1}}{\rho} \sin^{2} \biggr( \frac{\lambda l}{2}
\biggr) =0 \, ,}
\end{equation}
where $\frac{m_{1}}{\rho}\gg 1$. From the equations
\eqref{eq:determ-frq2}, \eqref{eq:determ-freq2} and
\eqref{eq:determ-freq-m}, we have the following properties of the
eigenvalues.

\begin{theorem}
\label{th:1}
If the discrete mass $m_{1}$ is located in the middle
$\frac{l}{2}$ of the integration interval $[0,l]$ of the
continuous-discrete boundary problem
\eqref{eq:motn-disc1}--\eqref{conjug-cond2-disc}, then all the even
eigenvalues of this problem coincide with the even
eigenvalues of the corresponding continuous boundary problem
\eqref{eq:motion-disc1}--\eqref{bound-cond-disc1}.
\end{theorem}

\begin{theorem}
\label{th:2}
For any location of the discrete mass $m_{1}$ at the integration
interval $[0,l]$ of the continuous-discrete boundary problem
\eqref{eq:motn-disc1}--\eqref{conjug-cond2-disc}, the first
eigenvalue of this problem is always smaller than the first
eigenvalue of the corresponding continuous boundary problem
\eqref{eq:motion-disc1}--\eqref{bound-cond-disc1}.
\end{theorem}

\begin{theorem}
\label{th:3} If the continuous-discrete system described by the
boundary problem \eqref{eq:motn-disc1}--\eqref{conjug-cond2-disc}
has only one discrete mass $m_{1}$ which is located in the middle
$\frac{l}{2}$ of the integration interval $[0,l]$, then the first
eigenvalue of the problem
\eqref{eq:motn-disc1}--\eqref{conjug-cond2-disc} is the smallest one
in comparison with the cases of any other location of the mass
$m_{1}$ at the integration interval.
\end{theorem}

Theorem~\ref{th:3} is in agreement with the physical sense of the
process under investigation: if the discrete mass $m_{1}$ is located
in the middle of the elastic body, then the amplitude of vibrations
on the first frequency at this point is maximal and, consequently,
the first eigenvalue and the first frequency of vibrations are
minimal.

In the case of two different discrete masses $m_{1}$ and $m_{2}$,
which are arbitrarily located at the cable, we have the next
equation to determine the eigenvalues of the problem
\eqref{eq:motn-disc1}--\eqref{conjug-cond2-disc}:

\begin{multline}
\label{eq:determ-freq-m2} \frac{1}{\lambda}\biggr[\sin(\lambda
l)-\frac{\lambda m_{1}}{\rho} \sin(\lambda b_{1})\sin(\lambda
(b_{1}+b_{2})) - \frac{\lambda m_{2}}{\rho}\sin(\lambda
b_{3})\sin(\lambda (b_{1}+b_{2})) \\
+ \frac{\lambda m_{1}m_{2}}{\rho}\sin(\lambda b_{1})\sin(\lambda
b_{2})\sin(\lambda b_{3})\biggr] =0 \, ,
\end{multline}
where $b_{i}$ is the length of the $i$-th interval of
continuity; $i=1, 2, 3$. At condition $b_{1}=b_{2}=b_{3}$
we can present equation \eqref{eq:determ-freq-m2} in a more simple form:

\begin{equation}
\label{eq:determ-freq-m3} \frac{1}{\lambda}\biggr[\sin(\lambda
l)-\frac{2\lambda (m_{1}+m_{2})}{\rho}
\sin^{2}\biggr(\frac{2l}{3}\biggr)\cos\biggr(\frac{\lambda l}{3}\biggr)
+ \frac{\lambda m_{1}m_{2}}{\rho}\sin^{3}\biggr(\frac{\lambda
l}{3}\biggr)\biggr] =0 \, .
\end{equation}

Finally, let us consider the general case of the continuous-discrete
boundary problem \eqref{eq:motn-disc1}--\eqref{conjug-cond2-disc},
when we have $n-1$ discrete masses $m_{i}$, $i=\overline{1,n-1}$, at
arbitrary points of the integration interval $[0,l]$. In this case
the determination of the eigenvalues is reduced to finding
of zeros of the function
\begin{equation}
\label{eq:determ-fnct} {\psi_{n+1}^{(1)}(\lambda)=
\psi_{n}^{(1)}S_{n}(b_{n})+\psi_{n}^{(2)}T_{n}(b_{n}) \, ,}
\end{equation}
with $\psi_{i+1}^{(1)}$ and $\psi_{i+1}^{(2)}$ given by the
recurrence formulas
\begin{equation}
\label{eq:recur1} {\psi_{i+1}^{(1)} =
\psi_{i}^{(1)}S_{i}(b_{i})+\psi_{i}^{(2)}T_{i}(b_{i}) \, ,}
\end{equation}
\begin{multline}
\label{eq:recur2} \psi_{i+1}^{(2)}=\frac{-\psi_{i}^{(1)}[\rho
\lambda^{2}T_{i}(b_{i})+m_{i}\lambda^{2}S_{i}(b_{i})]}{\rho} \\
+ \frac{\psi_{i}^{(1)}[\rho
S_{i}(b_{i})-m_{i}\lambda^{2}T_{i}(b_{i})]}{\rho} \, , \quad
i=\overline{1,n-1} \, ,
\end{multline}
where $S_{i}$ and $T_{i}$ are the normal fundamental functions
\begin{equation}
\label{fund-funct} {S_{i}(b_{i})=\cos (\lambda b_{i}) \, , \quad
T_{i}(b_{i})=\frac{\sin(\lambda b_{i})}{\lambda} \, ,}
\end{equation}
$\psi_{1}^{1}$ and $\psi_{1}^{2}$ are known from the left boundary
condition \eqref{bound-cond-1}; $b_{i}$ is the length of the $i$-th
interval of continuity, $i=\overline{1,n}$. As a result, we can
write the next equation for the determination of the eigenvalues:
\begin{equation}
\label{eqn-eigen} {\frac{1}{\lambda}[\sin(\lambda l)-\lambda
f_{n}(\lambda)]=0 \, ,}
\end{equation}
where $f_{n}(\lambda)$ includes all functions calculated by the
recurrence formulas \eqref{eq:recur1} and \eqref{eq:recur2}. From
\eqref{eqn-eigen} the following property of the eigenvalues can be
derived.

\begin{theorem}
The eigenvalues of the continuous-discrete boundary problem
\eqref{eq:motn-disc1}--\eqref{conjug-cond2-disc} are always less or
equal than the eigenvalues of the corresponding continuous boundary
problem \eqref{eq:motion-disc1}--\eqref{bound-cond-disc1}.
\end{theorem}

\subsection{The Problem for the Case When Only the Loads Move}

If only discrete loads with masses $m_{1}$, $m_{2}$, ...,
$m_{n-1}$ move along the cable-way at a steady speed, then the
mathematical model for such controlled dynamical system is the
next non-stationary boundary problem:

\begin{equation}
\label{eq:motion-disc} {\rho\frac{\partial^{2}u_{i}}{\partial
t^{2}}=T\frac{\partial^{2}u_{i}}{\partial x^{2}} \, ,}
\end{equation}

\begin{equation}
\label{bound-cond-disc} {u_{1}(0,t)=u_{n}(l,t)=0 \, ,}
\end{equation}
with conjugation conditions

\begin{equation}
\label{conjug-cond-disc1}
{u_{i+1}(l_{i},t)=u_{i}(l_{i},t)|_{x=l_{i}(t)}} \, ,
\end{equation}

\begin{equation}
\label{conjug-cond-disc2} m_{i}\left.\left(\frac{\partial^{2}u_{i}}{\partial
t^{2}}+2v\frac{\partial^{2}u_{i}}{\partial x\partial
t}+v^{2}\frac{\partial^{2}u_{i}}{\partial
x^{2}}\right)\right|_{x=l_{i}(t)}=\left.T\left(\frac{\partial u_{i+1}}{\partial
x}-\frac{\partial u_{i}}{\partial x}\right)\right|_{x=l_{i}(t)} \, ,
\end{equation}
where \, $l_{i}(t)=l_{i}+v t$, \, $l_{0}(t)=0$, \, $t_{0} \leq t
\leq t_{1}$.

Let us consider the boundary problem
\eqref{eq:motion-disc}--\eqref{conjug-cond-disc2} with only one moving load.
Using normal fundamental functions, we obtain the
following equation for the determination of non-stationary eigenvalues:
\begin{equation}
\label{eq:determ-freq} {\frac{\sin(\lambda
l)}{\lambda}-\frac{m_{1}}{\rho}(a^{2}+v^{2}) \sin (\lambda b_{1})
\sin (\lambda b_{2})=0 \, ,}
\end{equation}
where $\lambda = \lambda (t)$; \, $b_{1}=l_{1}$ and $b_{2}=l-l_{1}$
at a fixed moment of time. The zeros of the left-hand side of
\eqref{eq:determ-freq} are non-stationary eigenvalues at the time
$t_{0} \leq t \leq t_{1}$. It is possible at each fixed moment of
time to find non-stationary frequencies $\omega_{n}=\omega_{n}(t)$
at the same interval by the known formula $\omega_{k}=a
\lambda_{k}$, $k=1, 2, \ldots$

The role of the speed $v$ as a controlled parameter, which changes
the dynamical characteristics of the system, is evident from
\eqref{eq:determ-freq} and \eqref{eq:determ-freq2}. The presence of
the speed $v$ in \eqref{eq:determ-freq} implies:
\begin{proposition}
In the case of one moving load, the eigenvalues and the natural
frequencies of the dynamical system
\eqref{eq:motion-disc}--\eqref{conjug-cond-disc2} decrease with the
controlled parameter $v$, i.e. decrease when the velocity increases.
\end{proposition}

\subsection{The Mathematical Model for the Case of the Whole System Motion}

Now we assume that the whole system of the aerial cable-way
moves at some steady speed. The mathematical model for such hybrid
controlled system is a non-stationary boundary problem of second order PDEs.
It consists of the equation of motion \eqref{eq:motion} subject
to the boundary conditions \eqref{bound-cond}:
\begin{equation}
\label{eq:motion} \rho\left(\frac{\partial^{2}u_{i}}{\partial
t^{2}}+2v\frac{\partial^{2}u_{i}}{\partial x\partial t}
+ v^{2}\frac{\partial^{2}u_{i}}{\partial
x^{2}}\right)=T\frac{\partial^{2}u_{i}}{\partial x^{2}} \, ,
\end{equation}
\begin{equation}
\label{bound-cond} {u_{1}(0,t)=u_{n}(l,t)=0 \, ,}
\end{equation}
where $v$ is the steady speed of the system motion; $l_{i}$ is the
initial coordinate of the cable point with the fixed discrete load
$m_{i}$, $i=\overline{1,n-1}$; $l_{i}(t)=l_{i}+vt$ is the coordinate
at the moment $t$ of the cable point with the fixed discrete load
$m_{i}$, $i=\overline{1,n-1}$; $l_{0}(t)=l_{0}=0$;
$l_{n}(t)=l_{n}=l$; $u_{i}(x,t)$ are the functions of deviation from
equilibrium position for the point $x\in[l_{i-1}(t),l_{i}(t)]$ at
the moment $t$, $i=\overline{1,n}$.
Solutions of the problem \eqref{eq:motion}--\eqref{bound-cond}
must also satisfy the conjugation conditions
\eqref{conjug-cond1}--\eqref{conjug-cond2} at points $x=l_{i}(t)$:
\begin{equation}
\label{conjug-cond1} {u_{i+1}(l_{i},t)=u_{i}(l_{i},t)|_{x=l_{i}(t)}} \, ,
\end{equation}
\begin{equation}
\label{conjug-cond2} m_{i}\left.\left(\frac{\partial^{2}u_{i}}{\partial t^{2}}
+2v\frac{\partial^{2}u_{i}}{\partial x\partial t}
+v^{2}\frac{\partial^{2}u_{i}}{\partial
x^{2}}\right)\right|_{x=l_{i}(t)}=T\left.\left(\frac{\partial u_{i+1}}{\partial
x}-\frac{\partial u_{i}}{\partial x}\right)\right|_{x=l_{i}(t)} \, ,
\end{equation}
where $m_{1}$, $m_{2}$, \ldots, $m_{n-1}$ are the masses of the
discrete loads at arbitrary points $x=l_{i}(t)$,
$i=\overline{1;n-1}$.

Speed $v$, which may be changed within the limits of $0<v\leq
v_{1}$, is the controlled parameter for this real dynamical system.
The role of the speed $v$ may be shown for the problem
\eqref{eq:motion}--\eqref{conjug-cond2} even in the case of
absence of discrete loads.

If in the equation \eqref{eq:motion} we do not neglect the Coriolis
acceleration $\frac{\partial ^{2}u}{\partial x \partial t}$,
then the natural frequency is determined as follows:
\begin{equation}
\label{eq:determ-w} \omega_{k} = \frac{\pi k}{l}
\frac{a^{2}-v^{2}}{a} \, , \quad  a>v \, .
\end{equation}

In the case when we do not take into account the Coriolis
acceleration in this boundary problem, the natural frequencies are
given by
\begin{equation}
\label{eq:determ-w2} \omega_{k} = \frac{\pi k}{l}
\sqrt{\frac{a^{2}-v^{2}}{a}}
 \, , \quad  a>v \, .
\end{equation}

From \eqref{eq:determ-w}, \eqref{eq:determ-w2} one can see that the
speed $v$ may significantly decrease the natural frequencies. Here,
the obtained frequencies, in contrast to the problems considered
below, do not depend on time but depend on the speed only.

Besides, if we consider the case when the interval of integration
changes with a constant velocity, then the natural frequency for the
given time interval $t_{0} \leq t \leq t_{1}$ may be determined from
the formula
\begin{equation}
\label{eq:determ-w3} \omega_{k}(t) = \frac{\pi k}{l(t)}
  \frac{a^{2}-v^{2}}{a} \, , \quad a>v \, ,
\end{equation}
when the Coriolis acceleration is taken
into account, or by
\begin{equation}
\label{eq:determ-w4} \omega_{k}(t) = \frac{\pi k}{l(t)}
\sqrt{\frac{a^{2}-v^{2}}{a}}
 \, , \quad  a>v \, ,
\end{equation}
if the Coriolis acceleration is neglected.
From the formulas \eqref{eq:determ-w3} and \eqref{eq:determ-w4}
we get the following result.
\begin{theorem}
When the interval of integration changes with a constant velocity,
the speed $v$ of the cable motion decreases the frequencies of
vibrations for problem \eqref{eq:motion}--\eqref{conjug-cond2}.
Also, the decrease of the integration interval, under the condition
of invariance of the initial tension, always increase the
frequencies.
\end{theorem}


\section{Conclusions}

In this paper we consider several eigenvalue problems that model the
dynamics of some hybrid systems of aerial cable-ways with different
assumptions. We begin with an overview of the methods found in the
literature for similar eigenvalue problems. The method of the normal
fundamental systems turns out to be the most effective one for the
considered problems. This method is applied, and the dependence of
the dynamical characteristics of the systems on the controlled
parameter is established. However, many open problems still exist in
this field. First of all, the most general eigenvalue problem
\eqref{eq:motion}--\eqref{bound-cond} is not studied for the case of
presence of the discrete loads. Besides, for many applications it is
important to consider systems where the masses of the discrete loads
are decreasing with time \cite{Kach-teor}. All these problems are
supposed to be solved by the same method in future investigations.


\section*{Acknowledgements}

The first author is grateful to the partial financial support
provided by the Control Theory Group (\textsf{cotg}) of the Center
for Research in Optimization and Control (\textsf{CEOC}) of the
University of Aveiro, for participation in \emph{The International
Conference of Hybrid Systems and Applications}, held at the
University of Louisiana, Lafayette, LA, USA, May 22-26, 2006. The
hospitality and the good working conditions at the University of
Aveiro are also gratefully acknowledged. Also the authors are
grateful to Enrique H. Manfredini for the suggestions regarding
improvement of the text.


\end{document}